\documentclass[11pt]{amsart}
\usepackage{graphicx,a4wide}

\DeclareMathOperator{\per}{per}

\begin{document}

\author{Michael Baake} 
\address{Fakult\"{a}t f\"{u}r Mathematik,
  Universit\"{a}t Bielefeld,\newline \hspace*{\parindent}Postfach
  100131, 33501 Bielefeld, Germany}
\email{mbaake@math.uni-bielefeld.de}
\author{David Damanik}
\address{Department of Mathematics, Rice University, Houston, TX 77005, USA}
\email{damanik@rice.edu}
\author{Uwe Grimm}
\address{Department of Mathematics and Statistics,
  The Open University,\newline \hspace*{\parindent}Walton Hall, 
  Milton Keynes MK7 6AA, United Kingdom} 
\email{uwe.grimm@open.ac.uk}

\title{Aperiodic Order and Spectral Properties}

\maketitle

\section{Introduction}

The concept of \emph{order} is fundamental to human culture. It not only
underlies much of art and architecture, the scientific approach
to the understanding of our world is based on detecting and describing
order in Nature, in the form of laws of Nature described by mathematics.

Although humans have an instinctive understanding of order, it is
surprisingly difficult to define precisely what order is. An example
of order in Nature is provided by a perfect crystal, such as a
flawless diamond, in which atoms are ordered in a periodically
repeating pattern. That Nature can accommodate more complex forms of
order is well known, not the least since the discovery of
quasicrystals by Dan Shechtman \cite{Dan}, for which he received the
Wolf Prize in Physics in 2011 and eventually the Nobel Prize in
Chemistry in 2011, almost 30 years after his ground-breaking
discovery. The characteristic feature of these materials is the fact
that, while showing a similar degree of atomic order, they display
symmetries that are incompatible with a periodic arrangement of atoms.

The theory of \emph{aperiodic order} considers mathematical structures
that possess order without periodicity. While the advent of
quasicrystals has provided additional physical context to the
research, it dates back to the beginning of the twentieth century,
with the work of Harald Bohr on almost periodic functions \cite{B1,B2}.
It has since developed into a fascinating field of modern mathematics,
with links to many areas of mathematics such as dynamical systems, 
harmonic analysis, spectral theory, number theory, to name but a few.
For a gentle introduction to the field, we refer to \cite{BGM};
a more comprehensive introductory account is given in \cite{TAO}.

Apart from the visual attraction of aperiodic tilings, an attractive
aspect of the field is the fact that one can make seemingly simple
statements which are easy to understand but turn out to be very
difficult to prove (similar to, say, Goldbach's conjecture in number
theory, which states that every even number can be written as the sum
of two prime numbers).  An example of such a question is whether there
exists a planar shape that can tile the entire plane without gaps or
overlaps (like in a puzzle), but does not admit any periodic
tilings. The answer to this question is still open, although there has
been some recent progress towards an answer by the Australian
mathematician Joan Taylor; see \cite{ST} as well as
\cite[Ex.~6.6]{TAO} and references therein.

In this exposition, we introduce the general idea behind aperiodic order
by means of simple but instructive examples, and provide a hint of why
spectral properties are of interest in this context. In doing so, we
will gloss over any technical details; we refer to \cite{TAO} for a
proper mathematical treatment.

\section{Point sets}

We introduce the notion of periodic and aperiodic order by considering
simple examples of point sets on the line, that is, in one dimension
of space. Sometimes we would like to distinguish different types of
points, say by assigning a colour to each point; each point is then
characterised by its position $x\in\mathbb{R}$ on the real line and by
its colour.  We start with a simple case where all points are located
at \emph{integer} positions along the line.  Imagine placing a red
point at any integer position of the real line to obtain the point set
$P_{0}$ of red points at all positions $n\in\mathbb{Z}$, which looks
like\medskip
\[
   \includegraphics{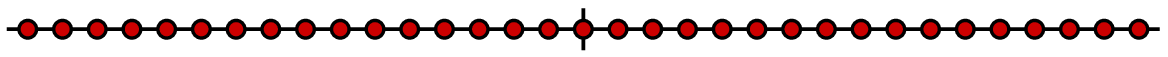}
\]
where the small vertical line denotes position $n=0$.  This point set,
which you have to imagine to continue indefinitely in both directions,
is \emph{periodic} with period $1$, because shifting all positions by
$1$ reproduces the same point set, so $P_{0}+1=P_{0}$ (where the
notation means adding $1$ to the position of each point in the
coloured point set $P_{0}$).  Of course, if shifting by $1$ maps the
point set onto itself, so does shifting by $2$, or indeed by
\emph{any} integer $n\in\mathbb{Z}$, so the set of periods of $P_{0}$
is $\per(P_{0})=\mathbb{Z}$. This reasoning holds true for \emph{any}
periodic point set; so once a point set possesses a period, it
automatically possesses an infinite set of periods, which forms a
\emph{lattice}, consisting of all integer linear combinations of a set
of fundamental periods. Here, we are in one dimension, and the
fundamental period (which is the smallest non-trivial period, where
$0$ is deemed to be a trivial period) is $1$.

Now, let us take every point at a  position $n$ with $n\equiv 1\bmod 4$
(which means that dividing $n$ by $4$ gives a remainder of $1$) and 
change its colour to blue, and call the corresponding point set $P_{1}$. 
It looks like this,\medskip
\[
   \includegraphics{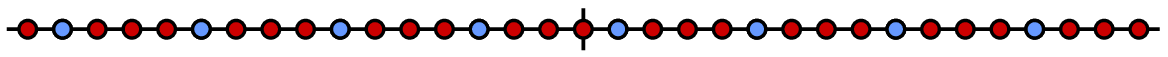}
\]
where again you have to imagine the point set to continue in both
directions, so we have changed the colour of infinitely many
points. What is the periodicity of the new point set? We now need to
shift $P_{1}$ by multiples of $4$ to respect the positions and the
colourings of points, so $\per(P_{1})=4\mathbb{Z}$.

As the next step, let us look at all points at positions $n$ with
$n\equiv 7\bmod 16$. As you can easily convince yourself, all these
points are currently red, so let us change them to blue to obtain the next
point set, which we call $P_{2}$. The result is\medskip
\[
   \includegraphics{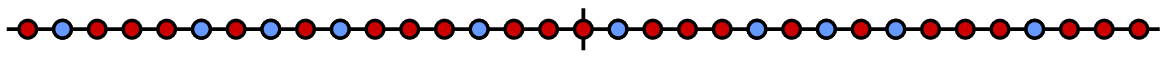}
\]
and is still periodic, but now only under shifts by multiples of $16$, so
$\per(P_{2})=16\mathbb{Z}$. 

We can continue this game, for example by defining $P_{k+1}$ as the
point set obtained from $P_{k}$ by changing the colour of all points
at positions $n\equiv (2\cdot 4^{k}-1)\bmod 4^{k+1}$. All these points
are still red in $P_{k}$, so will then become blue. The resulting
point set $P_{k+1}$ is then periodic under shifts by multiples of
$4^{k+1}$. In this way, we obtain point sets $P_{k}$ for all
integer $k\ge 0$, which are periodic with
$\per(P_{k})=4^{k}\mathbb{Z}$. With increasing $k$, the periods get
sparser and sparser (indeed, in order to actually see this happening
even for the next step $k=3$ you would need to consider a longer part
of the point set than we displayed above, because this step affects
points at positions $n\equiv 31\bmod 64$ only), and if we kept on
performing this process indefinitely and consider the limiting case
where $k$ becomes infinite, we eventually end up with a point set $P$
that no longer has any periods at all.

Such a point set that does not admit \emph{any} non-trivial period is
called \emph{non-periodic}. In fact, our point set $P$ is not just
non-periodic, but actually \emph{aperiodic} in the sense of
\cite{TAO}, which is a somewhat stronger statement. A proper
explanation of the difference between the notions of non-periodicity
versus aperiodicity requires a more careful definition of limits of
shifted point sets, but the main point of the distinction is to
eliminate certain `trivial' situations. An example is the case where
you take our original single-colour point set $P_{0}$ and change the
colour of a single point, say the point at $0$. The resulting point
set is non-periodic, because any non-zero shift moved the blue point
at $0$ to another position and hence changes the point set, because
there is no other blue point in the entire point set. However, if you
keep shifting the blue point further and further away, the point set
will look more and more like the original point set $P_{0}$ around the
origin, and in the limit where the blue point has been moved off to
infinity, the periodicity of the point set is restored. This point set
would not be considered as aperiodic, because periodicity is violated
only locally, not globally. Because, in our construction of the point
sets $P_{k}$, we change the colour of infinitely many points in
each step, this is not the case here, and $P$ possesses the stronger
aperiodicity property.

Although the point set $P$ is aperiodic, it is clearly ordered in some
sense, because it is built from an explicit construction which
determines the colour for each position uniquely. Even if you did not
know where the origin was located, you can still recognise this
order. For instance, if you pick a red point which is located between
two blue points, you know immediately that every second point along
either direction will be red as well, because all points at even
positions stay red in our construction. This also shows that you can
never find two blue points next to each other. However, like for the
simple periodic set $P_{0}$, if you do not know where the origin is
located, you actually cannot decide where the origin would have been
from looking at an arbitrarily large finite part of the set, because
any local arrangements of colours repeat indefinitely often along the
line, just not periodically. This property is called
\emph{repetitivity} and arises here as a consequence of the systematic
way we used to perform the colour changes, which affected points in
the same way anywhere along the line.

So our point set $P$ is an example of a structure that is both ordered
and aperiodic. It is closely related to a class of sequences known as
Toeplitz sequences \cite{T}. The theory of \emph{aperiodic order} is
concerned with understanding such point sets and analysing their
properties.

\section{Substitution and inflation}

You may wonder why we used the rather specific way of changing colours
in our construction of the point sets $P_{k}$ above. Clearly, there
are lots of ways to produce aperiodic point sets in a similar way.
The reason why we chose this particular approach is that the point
set $P$ is a nice and simple example of a set that can also be
obtained in a different way, namely by what is known as a
\emph{substitution} or \emph{inflation} rule (where the latter is
commonly used for the geometric interpretation which will be discussed
later). To see this, let us denote the sequence of the two coloured
points by letters $r$ (for red) and $b$ for blue, and consider the
rule $S$ that maps $r\mapsto rb$ and $b\mapsto rr$.  Applying
this rule repeatedly, starting from a single letter, gives
\[
   r\,\stackrel{\scriptstyle S}{\longmapsto}\, 
   rb\,\stackrel{\scriptstyle S}{\longmapsto}\, 
   rbrr \,\stackrel{\scriptstyle S}{\longmapsto}\, 
   rbrrrbrb \,\stackrel{\scriptstyle S}{\longmapsto}\, 
   rbrrrbrbrbrrrbrr \,\stackrel{\scriptstyle S}{\longmapsto}\, \dots 
\]
where, in each step, every letter is replaced by a pair of letters
according to the rule $S$. You can continue to do this as long as you
like, producing longer and longer words in the two letters $r$ and
$b$, and in the limit you obtain an infinite word, $v$ say, that is
mapped onto itself under the rule, so $Sv=v$. This word is thus
invariant under the application of the rule $S$, and in this sense
possesses a symmetry under this operation, sometimes referred to as an
\emph{inflation symmetry}. The surprising result is that the infinite
sequence $v$ you obtain in this way will \emph{exactly} reproduce the
sequence of colours in $P$, starting with the red point at position
$n=0$. To mathematically prove why this is so requires some work; if
you are interested, you can find the argument in
\cite[Chapter~4.5.1]{TAO}, where this example is referred to as the
\emph{period doubling substitution}.

You may ask what happens to the `left' part of the point set $P$,
consisting of all points at negative positions $n<0$. In fact, you can
obtain this in a very similar fashion, starting from a two-letter seed and
keeping track of the origin as follows
\begin{align*}
 r|r\,  & \stackrel{\scriptstyle S}{\longmapsto}\, rb|rb 
 \stackrel{\scriptstyle S}{\longmapsto}\, 
 rbrr|rbrr \,\stackrel{\scriptstyle S}{\longmapsto}\,
 rbrrrbrb|rbrrrbrb\\ 
 &\stackrel{\scriptstyle S}{\longmapsto}\, rbrrrbrbrbrrrbrr|rbrrrbrbrbrrrbrr 
 \,\stackrel{\scriptstyle S}{\longmapsto}\, \dots 
\end{align*}
which you can again continue indefinitely. Considering every second
word, starting from the seed $r|r$, you find that the words grow into
a bi-infinite word $w$ which coincides with $v$ on the `right' of the
origin, and which satisfies $S^{2}w=w$. It produces \emph{precisely}
the sequence of colours of the point set $P$. Note that, if you had
chosen the other word obtained by using an odd number of applications
of $S$ on $r|r$, the only difference would be that the point at
position $n=-1$ is blue rather than red, all other colours remain the
same.
 
Substitution rules like $S$ have been studied extensively, and produce
many well-known examples of interesting sequences. The most famous
such sequence is named after Leonardo of Pisa, also known as
Fibonacci, who introduced it in his book \emph{Liber Abaci} already in
1202, although it was apparently familiar to Indian mathematicians
even earlier. It was motivated by studying the evolution of a rabbit
population, with the rule that, in one step, any adult rabbit produces
one offspring, and any juvenile rabbit matures to an adult
rabbit. This is, of course, a very simplified model in which rabbits
live and (asexually) reproduce eternally, and the total population
grows exponentially! Let us denote the adult rabbits by $\ell$ (for
large) and the young rabbits by $s$ (for small), the Fibonacci rule
$F$ is $\ell\mapsto\ell s$ and $s\mapsto\ell$. Applying the rule
repeatedly, starting with a single adult, gives
\[
\ell \, \stackrel{\scriptstyle F}{\longmapsto}\, 
\ell s \,  \stackrel{\scriptstyle F}{\longmapsto}\, 
\ell s \ell \, \stackrel{\scriptstyle F}{\longmapsto}\, 
\ell s \ell \ell s \, \stackrel{\scriptstyle F}{\longmapsto}\, 
\ell s \ell \ell s \ell s \ell \, \stackrel{\scriptstyle F}{\longmapsto}\, 
\ell s \ell \ell s \ell s \ell \ell s \ell \ell s \,
\stackrel{\scriptstyle S}{\longmapsto}\, \dots
\]
which, when repeating the process indefinitely, produces an infinite
word $v$ which satisfies $Fv=v$ and is known as the \emph{Fibonacci
  sequence}.

Recognising that each finite word in the iteration above is the
concatenation of the two previous words, the number of letters of any
one of these words is the sum of the number of letters of the two
previous words. This produces the sequence of \emph{Fibonacci numbers}
$1,2,3,5,8,13,21,34,\ldots$ satisfying the recursion relation
$f_{k+1}=f_{k}+f_{k-1}$, with initial conditions $f_{0}=0$ and
$f_{1}=1$ (in which case the list above starts with $f_{2}=1$, and
$f_{3}=f_{2}+f_{1}=1+1=2$ and so on). The Fibonacci numbers thus give
the total number of rabbits after a number of generations.  Counting
the numbers of adult or young rabbits, so either of the letter $\ell$
or the letter $s$, in each of generations again produces the same
sequence, so that in a word of length $f_{k+1}$ there are exactly
$f_{k}$ letters $\ell$ (adult rabbits) and $f_{k-1}$ letters $s$
(young rabbits). Using this observation, it is not difficult to show
that the ratio of letters $\ell$ and $s$ (the ratio of adult to young
rabbits), as the number of generations grows, approaches the limit
\[
    \lim_{k\to\infty}\frac{f_{k}}{f_{k-1}} \, = \, \frac{1+\sqrt{5}}{2}
    \, = \, 1.6180339887\ldots
\]  
which is an irrational number, conventionally denoted by $\tau$, known
as the \emph{golden ratio} (and an important number in art and
architecture, representing an `ideal' way of dissecting an interval
into two parts).  The fact that $\tau$ is irrational shows that the
Fibonacci sequence $v$ cannot be periodic. Indeed, assuming that $v$
would repeat periodically after, say, $N$ letters, the ratio of
letters in $v$ would have to be the same as their ratio in a finite
word of length $N$, and hence a rational number with a denominator of
at most $N$.

As for the rule $S$, we can also produce a two-sided sequence by
taking every second step in the iteration of $F$ on the two-letter
seed $\ell|\ell$, giving
\[
   \ell|\ell \, \stackrel{\scriptstyle F^{2}}{\longmapsto}\,
   \ell s \ell | \ell s \ell  \, \stackrel{\scriptstyle F^{2}}{\longmapsto}\,
   \ell s \ell \ell s \ell s \ell | \ell s \ell \ell s \ell s \ell 
   \, \stackrel{\scriptstyle F^{2}}{\longmapsto}\, \dots
\]
which produces a bi-infinite word $w$ which coincides with $v$ on the
`right' and satisfies $F^{2}w=w$.

There is a natural way to interpret the Fibonacci sequence as a point
set on the real line, in a way that the rule $F$ becomes an
\emph{inflation} rule in the following sense. Let us associate to the
two letters $\ell$ and $s$ two interval lengths, a long one (to
$\ell$, to fit the adult rabbits in) and a short one (to $s$, for the
young rabbits). A natural way to choose the length is given by the
golden ratio again, so let us choose the length of the interval $\ell$
to be $\tau$ and the length of the interval $s$ to be $1$ (for the
mathematical reason for this choice see 
the discussion of geometric inflation rules in \cite[Ch.~4]{TAO}).
Then, the geometric interpretation of the rule $F$ is\medskip
\[
   \includegraphics[width=0.85\textwidth]{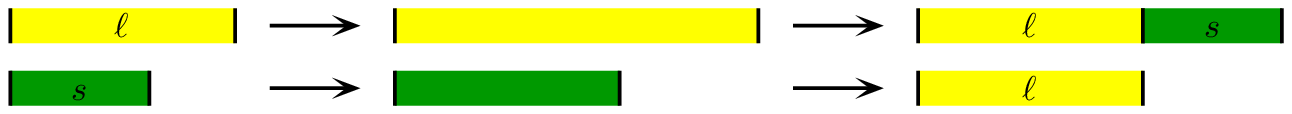}
\]
which consists of a scaling of the intervals by a factor $\tau$,
followed by the dissection of the long interval into a long and short
one (according to the rule $\ell\mapsto\ell s$, which is geometrically
consistent because $\tau^2=\tau+1$) and interpreting the scaled short
interval as a long one (according to the rule $s\mapsto \ell$). The
geometric versions of the infinite or bi-infinite words $v$ and $w$
become series of intervals\medskip
\[
   \includegraphics[width=0.85\textwidth]{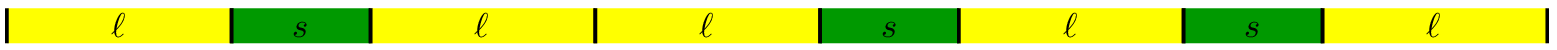}
\]
which are invariant under this geometric inflation map.

\section{Cut and project sets}

\begin{figure}[b]
\includegraphics[width=0.85\textwidth]{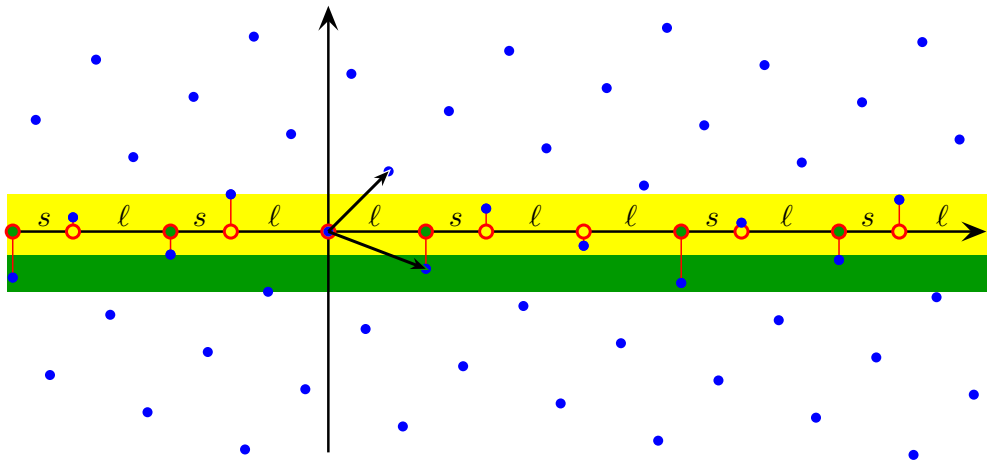}
\caption{Cut and project description of the Fibonacci sequence from a
  planar lattice (blue dots). Lattice points within the yellow strip
  are projected onto left endpoints of long intervals $\ell$, while
  points within the green strip project to left endpoints of short
  intervals $s$. The resulting sequence of long and short intervals
  gives the Fibonacci sequence.\label{fig:fiboproj}}
\end{figure}

The Fibonacci sequence, in its geometric presentation, can also be
obtained in a seemingly very different way, which is sketched in
Figure~\ref{fig:fiboproj}. Here, the long and short intervals are
obtained by a projection of a two-dimensional periodic point lattice
(the blue dots) onto a one-dimensional subspace, selecting all the
points that fall within the coloured strip. Note that points within
the yellow part of the strip give rise to left endpoints of long
intervals, while those within the green part correspond to left
endpoints of short intervals. With the chosen setup, in particular the
location and width of the strip used for the selection of lattice
points, the projected one-dimensional tiling turns out to be exactly
the same as the one obtained from the inflation description discussed
above.  Note that the choice of lattice is not unique; our choice is
motivated by an interesting connection to number theory, namely the
Minkowski embedding of the ring $\mathbb{Z}[\tau]=\{m+n\tau\mid
m,n\in,\mathbb{Z}\}$; see \cite[Ch.~3.4]{TAO} for details.

This interpretation of the Fibonacci case in terms of a \emph{cut and
  project} set (or \emph{model set}, a notion that goes back to
Meyer's work on harmonious sets \cite{Meyer}) is another indication of
the inherent order that is `hidden' in the aperiodic
sequence. Although it is not periodic, it is very closely related to a
periodic structure, albeit in two rather than in one dimension. This
is an important property, and the cut and project construction can be
generalised and applied in a quite general setting. The resulting cut
and project sets are now quite well understood, because the underlying
higher-dimensional periodicity provides a \emph{quasiperiodic} order,
which is a particular case of the general notion of almost periodic
order.

The fact that the Fibonacci sequence allows both an inflation and a
projection description should not mislead you to assume that this
happens in general. Indeed, we are looking at a very special situation
here, although many of the `nice' examples will be of this kind. Given
an inflation rule, it does not automatically allow for an embedding
into a periodic lattice in a higher-dimensional space; only certain
inflation rules can be represented in this way at all. It turns out
that the period doubling sequence discussed at the beginning does in
fact have a projection description, but only in a setting where the
periodic lattice lies in a more general space, which is not a
finite-dimensional Euclidean space, but in general even this is not
guaranteed. Conversely, given a cut and project description, the
projected structure does not automatically possess an inflation
symmetry. In the Fibonacci setup shown in Figure~\ref{fig:fiboproj},
this is only true if the strip is chosen appropriately.

\section{Spectral properties}

The one-dimensional examples discussed above should provide an
intuitive idea about the type of structures that we have in mind when
we talk about aperiodic order. The question then arises of how one can
characterise order in such aperiodic structures. This is where
spectral properties find an application, inspired by applications in
crystallography, physics and materials science.

Crystalline order is defined in terms of the \emph{diffraction} of a
material. Diffraction is experimentally observed as the pattern of
radiation (such as X-rays) that is scattered by the material. It
mainly measures two-point correlations in a structure, with a
point-like diffraction pattern indicating strong (crystalline) order,
and diffuse scattering typically being interpreted by an absence of
order.  Mathematically, the quantity that is probed (in the simplest
approximation of what happens in reality) is the
\emph{autocorrelation}. Let us consider an example. Take our point set
$P$ from our first example, and interpret all red points as scatterers
(for simplicity, we associate a scattering strength $u(r)=1$ to these)
and consider the blue points as empty (scattering strength
$u(b)=0$). The autocorrelation for a given distance $m$ is then the
average over the product of scattering strengths of points at distance
$m\in\mathbb{Z}$, so
\[
     a(m) \, = \, \lim_{N\to\infty}\frac{1}{2N+1}\sum_{n=-N}^{N}
     u(w_{n}) u(w_{n+m})
\]
where $w_{n}$ denotes the letter ($r$ or $b$) at position $n$ in the
bi-infinite sequence $W$. So $0\le a(m)\le 1$ is the average proportion of
times you find two scatterers at distance $m$ along the line.

\begin{figure}
\includegraphics[width=0.85\textwidth]{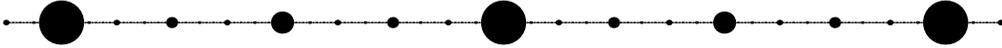}
\caption{Sketch of the diffraction pattern associated to 
the point set $P$. Here, contributions at $k$ are represented 
by a disk, centred at $k$, of an area that is proportional to 
the diffraction intensity $I(k)$.\label{fig:pddiff}}
\end{figure}

The autocorrelation coefficients $a(m)$ can be explicitly calculated
for this case, giving $a(0)=2/3$ and
\[
    a(m) \, = \, \frac{2}{3} \left(1-\frac{1}{2^{r+1}}\right)
\]
for $m=(2\ell+1)2^{r}$ with $r\ge 0$ and $\ell\in\mathbb{Z}$.  The
\emph{diffraction} is obtained from the autocorrelation by what is
called a Fourier transformation, which essentially is analysing the
frequency distribution of the autocorrelation \cite{H}. This is
similar to the frequency analysis of a sound. For the example at hand,
the increasingly long periods of powers of $4$ in the point set $P$
give rise to frequencies at rational numbers with powers of $2$ in the
denominator. The diffraction in this case is thus concentrated on a
dense point set consisting of all integers $k\in\mathbb{Z}$, with
intensity $I(k)=4/9$, as well as all rational numbers $k=(2n+1)/2^{r}$
with $n\in \mathbb{Z}$ and $r\ge 1$ (whose denominators are powers of
$2$), with the diffraction intensity at $k$ given by
\[
   I(k) \, = \, \frac{1}{9\cdot 4^{r-1}}\, ;
\]
see \cite[Ch.~9.4.4]{TAO} for details.  So, despite the fact that the
positions are dense in space, the decaying nature of $I(k)$ as a
function of $r$ means that contributions above any given threshold
involve a discrete set of positions only, and in any experiment you
would see a discrete patterns of spots of different intensity. A
sketch of the diffraction pattern, which is periodic with period $1$,
is shown in Figure~\ref{fig:pddiff}. As mentioned previously, the
point set $P$ admits an interpretation as a cut and project set, so
the pure point nature of the diffraction spectrum is in line with the
general result \cite{Sch} that such sets, under rather general
assumptions, are pure point diffractive.

Diffraction is one spectral measure attached to our point set; we
shall briefly mention two other spectral measures that have been
studied extensively. The first of these is the \emph{dynamical
  spectrum}. This is related to the action of translations on our
point set, and considering the space of all point sets obtained by
such translations and appropriate limits of translates.  The idea of
associating a spectral measure to this dynamical system goes back to
Koopman \cite{K} and von Neumann \cite{vN}. For a detailed discussion of
dynamical spectra for substitution-based systems we refer to
\cite{Q}. It turns out that the dynamical spectrum is generally richer
than the diffraction spectrum, which is due to the fact that it can
detect order beyond the two-point correlations that diffraction can
see. It has been known for a long time that there is a close
connection between diffraction and dynamical spectra; indeed, the
original proof that model sets are, under rather general assumptions,
pure point diffractive \cite{Sch} was based on an argument
linking the two spectra in the pure point case.  Recent work
\cite{BLvE} has further elucidated the connection between the dynamical
and diffraction spectra, which now is reasonably well understood.
 
Finally, another spectral measure motivated by physics is given by the
energy spectra of Schr\"{o}dinger operators describing the motion of
electrons in a solid. In a periodic system, valence electrons can move
`freely', in what is known as a Bloch wave, whereas in a disordered
system, where certain positions are energetically favourable to
others, electrons would be expected to be localised at such
places. Aperiodically ordered structures are somewhere in between
these two, because on the one hand they have motives that keep
repeating throughout the systems, while on the other hand they lack
periodicity which could sustain Bloch wave solutions. Indeed, it turns
out that for large classes of aperiodic models, one finds a very
peculiar behaviour, in which electrons are neither localised nor can
move freely, and associated to this behaviour are particularly
interesting spectral properties. We refer to \cite{DEG} 
and references therein for a recent comprehensive review of
the results in this area, and to \cite{DGY} for a detailed
analysis of the Fibonacci case. 

Although quite a bit is now known about the spectral properties of
Schr\"{o}dinger operators for large classes of one-dimensional
examples, there is currently no satisfactory understanding of the
relation between the spectral properties of these systems and the
other two spectral properties discussed earlier, if indeed such a
relation exists. In some respect, these spectral measures behave in
opposite ways; in diffraction and dynamical spectra, point spectra
indicate order in a system, while in the Schr\"{o}dinger case a
periodically ordered system shows a continuous spectrum. It is our
hope that further investigation of aperiodically ordered systems may
shed some light on this open question.

\begin{figure}
\includegraphics[width=0.85\textwidth]{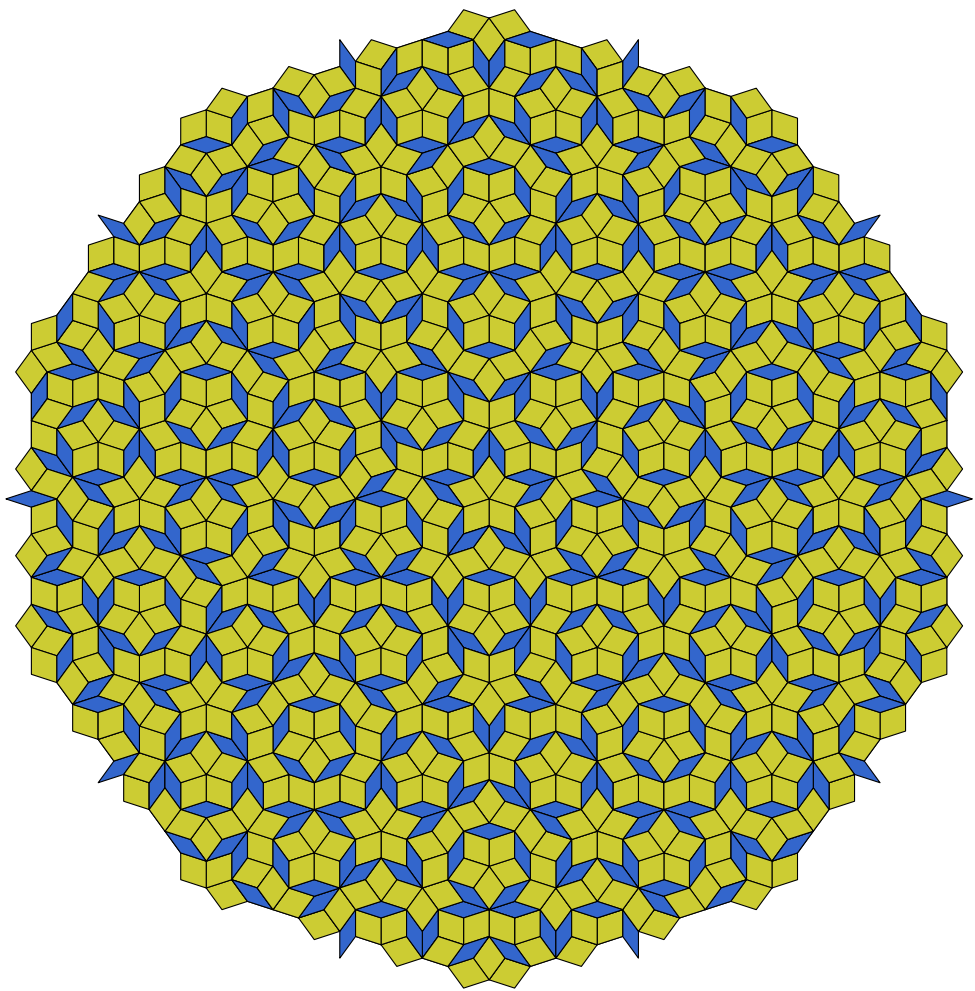}
\caption{A patch of the rhombic Penrose tiling.\label{fig:pen}}
\end{figure}

\section{Summary}

\begin{figure}
\includegraphics[width=0.85\textwidth]{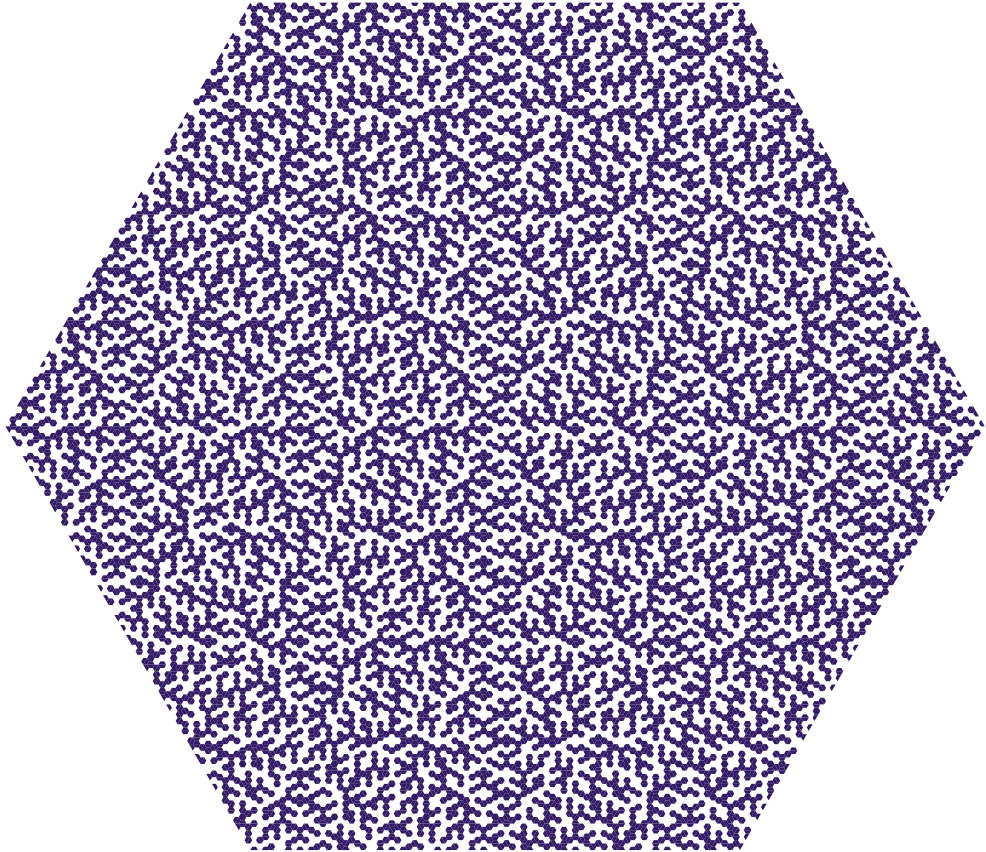}
\caption{A patch of Joan Taylor's llama tiling.\label{fig:llama}}
\end{figure}

In this snapshot, we limited the discussion to one-dimensional point
sets (or tilings by intervals).  However, aperiodic order is not
limited to one dimension, quite on the contrary some of the beauty of
the subject becomes apparent when considering aperiodic tilings in
higher dimensions.  The most famous example is Penrose's tiling
\cite{P}, of which there exist a number of variants; see
\cite[Ch.~6.2]{TAO}. A patch of the version with two rhombic tiles is
shown in Figure~\ref{fig:pen}.  Similar to the one-dimensional
Fibonacci system discussed above, Penrose's tiling can be described
either as a two-dimensional inflation tiling or as a cut and project
set, in this case using a lattice in at least four-dimensional space.
Its diffraction is pure point and shows perfect tenfold symmetry,
which is incompatible with periodicity by the crystallographic
restriction \cite[Ch.~3.2]{TAO}.

A more recent, stunning example is Joan Taylor's llama tiling shown in
Figure~\ref{fig:llama}, which is related to the open question
mentioned at the beginning. Here, the tiles are hexagons of two
different types (one of which is kept white in the figure), and the
name refers to the fact that the smallest connected cluster of tiles
of one colour form little a shape that is resembles the outline of a
llama. For more details about this tiling, we refer to
\cite[Ex.~6.6]{TAO}.

\bigskip

\end{document}